\documentclass[11pt]{amsart}

\usepackage{amsmath}
\usepackage{amsxtra}
\usepackage{amscd}
\usepackage{amsthm}
\usepackage{amsfonts}
\usepackage{amssymb}
\usepackage{eucal}
\usepackage[all]{xy}
\usepackage{graphicx}
\usepackage{tikz-cd}
\usepackage{mathrsfs}
\usepackage{euler}
\usepackage{hyperref}
\usepackage{color}
\usepackage{longtable}
\usepackage{float}
\usepackage{caption}
\usepackage{enumerate}
\usepackage{faktor}
\usepackage{lineno}
\usepackage{tikz}
\usepackage{tikz-cd}
\usetikzlibrary{matrix,arrows,decorations.pathmorphing}
\usepackage{kantlipsum}

\theoremstyle{plain}
  \newtheorem{theorem}{Theorem}[section]

  \newtheorem{lemma}[theorem]{Lemma}
  
  \newtheorem{corollary}[theorem]{Corollary}

\theoremstyle{definition}
	\newtheorem{definition}[theorem]{Definition}

	\newtheorem{non-example}[theorem]{Non-example}
	
	\newtheorem{conjecture}[theorem]{Conjecture}
	
\theoremstyle{remark}
	\newtheorem{remark}[theorem]{Remark}

\numberwithin{equation}{section}

\DeclareMathAlphabet{\mathcal}{OMS}{cmsy}{m}{n}

\newcommand{\N}{\mathbb{N}}
\newcommand{\M}{\mathbb{M}}

\renewcommand{\P}{\mathbb{P}}
\newcommand{\Q}{\mathbb{Q}}
\newcommand{\R}{\mathbb{R}}

\newcommand{\sca}{\mathscr{A}}

\newcommand{\scm}{\mathscr{M}}
\newcommand{\sco}{\mathscr{O}}

\newcommand{\scx}{\mathscr{X}}

\newcommand{\calh}{\mathcal{H}}

\newcommand{\calz}{\mathcal{Z}}
\newcommand{\mfp}{\mathfrak{p}}

\newcommand{\mfR}{\mathfrak{R}}

\newcommand{\spec}{\mathrm{Spec}\,}

\newcommand{\supp}{\mathrm{Supp}}

\newcommand{\gal}{\mathrm{Gal}}

\newcommand{\rad}{\mathrm{rad}}

\newcommand{\pgl}{\mathrm{PGL}}

\newcommand{\na}{\mathrm{na}}
\newcommand{\rat}{\mathrm{Rat}}
\newcommand{\crit}{\mathrm{crit}}
\newcommand{\Crit}{\mathrm{Crit}}

\newcommand{\preper}{PrePer}

\newcommand{\paren}[1]{\left(#1\right)}
\newcommand{\set}[1]{\left\{#1\right\}}
\newcommand{\sbrac}[1]{\left[#1\right]}
\newcommand{\abrac}[1]{\left<#1\right>}
\newcommand{\verts}[1]{\left\lvert#1\right\rvert}

\newcommand\restr[2]{{
  \left.\kern-\nulldelimiterspace 
  #1 
  \right|_{#2} 
  }}

\title{The $abcd$ conjecture, uniform boundedness, and dynamical systems}

\author{Robin Zhang}
\address{Department of Mathematics, Columbia University \newline
	\indent Department of Mathematics, Massachusetts Institute of Technology}
\email{rzhang@math.columbia.edu, robinz@mit.edu}

\date{March 6, 2024}

\begin{document}

\begin{abstract}
	We survey Vojta's higher-dimensional generalizations of the $abc$
	conjecture and Szpiro's conjecture as well as
	recent developments that apply
	them to various problems in arithmetic dynamics. In particular,
	the ``$abcd$ conjecture'' implies a dynamical
	analogue of a conjecture on the uniform boundedness of torsion
	points and a dynamical analogue of Lang's conjecture on lower
	bounds for canonical heights.
\end{abstract}

\maketitle

\dedicatory{To the memory of Lucien Szpiro (1941--2020)}


\section{Szpiro to \texorpdfstring{$abc$}{abc}}

The $abc$ conjecture, which originated in a conversation between Masser and
Oesterl\'{e} in 1985 as an approach to Szpiro's conjecture
(cf. \cite[Section~3]{oesterle}), has
been popularly described as ``the most important unsolved problem in
Diophantine analysis''~\cite{goldfeld1996}. Many articles have
been written about their implications across
number theory and Diophantine geometry
(e.g. the wonderful surveys
\cite{goldfeld1996, goldfeld1999, granville-tucker,
waldschmidt}).

This article highlights recent developments that extend this
discussion, with applications of the ``$abcd$ conjecture''
to various problems in arithmetic dynamics.
This $abcd$ conjecture lies between Vojta's higher-dimensional
``$abcde\ldots$ conjecture'' and the classical $abc$ conjecture.
It has been shown by Looper~\cite{looper1, looper2} to imply the
uniform boundedness of preperiodic points of polynomials on $\P^1$
and a weak version of the dynamical Lang conjecture on
points of small canonical height. The conditional uniform boundedness
of preperiodic points also yields the conditional non-existence of
rational (or even quadratic) periodic points of unicritical
polynomials of large enough degree.

\subsection{Szpiro's conjecture}
Let $E$ be an elliptic curve over the rational numbers $\Q$ with
global minimal Weierstrass equation
\[
	E: y^2 + a_1 xy + a_3 y = x^3 + a_2 x^2 + a_4 x + a_6.
\]

\noindent
Tate's notes~\cite{tate} (based on his letter to Cassels) define
invariants of $E$,
\begin{align*}
	c_2 &= a_1^2 + 4a_2, \\
	c_4 &= a_1 a_3 + 2a_4, \\
	c_6 &= a_3^2 + 4a_6, \\
	c_8 &= a_1^2 a_6 - a_1 a_3 a_4 + 4a_2 a_6 + a_2 a_3^2 - a_4^2.
\end{align*}
The discriminant $\Delta_E$ is $-c_2^2 c_8 - 8c_4^3 -27c_6^2
+ 9c_2 c_4 c_6$. This discriminant is minimal by our minimality
assumption on the model.
The conductor $N_E$ of $E$ is the product
\[
	N_E := \prod_{p \textrm{ prime}} p^{\epsilon_p + \delta_p}.
\]
The tame part $\epsilon_p$ of the conductor is
\[
	\epsilon_p :=
		\begin{cases}
			0 & \textrm{if } E \textrm{ has good reduction at } p, \\
			1 & \textrm{if } E \textrm{ has multiplicative reduction at } p, \\
			2 & \textrm{if } E \textrm{ has additive reduction at } p.
		\end{cases}
\]
The wild part $\delta_p$
comes from the $\ell$-adic Swan representation
and is nonzero only if $p < 5$
(for more about the conductor of elliptic curves,
see Silverman~\cite[Section~IV.10]{silverman-advanced}).

Following Szpiro's work on the Shafarevich and Mordell conjectures
over function fields of positive characteristic,
he hoped to use Arakelov theory to prove the effective Mordell
conjecture (see \cite[p.~1771--1774]{notices-szpiro}).
Szpiro's conjecture was presented at a talk in Hannover in 1982 in
the following form.
\begin{conjecture}[{Weak Szpiro's conjecture \cite[Conjecture~1]{oesterle}}]
	\label{conj:szpiro-weak}
	There exist positive real numbers $\alpha$ and $\beta$
	such that
	\[
		\verts{\Delta_E} \leq \alpha N_E^\beta,
	\]
	for any elliptic curve $E$ over $\Q$ with minimal discriminant
	$\Delta_E$ and conductor $N_E$.
\end{conjecture}

With a lower bound on $\beta$ and with an expression of $\alpha$ as a
function of $\beta$, there are stronger versions of Szpiro's
conjecture.
\begin{conjecture}[{Strong Szpiro's conjecture \cite[Conjecture~2]{oesterle}}]
	\label{conj:szpiro-strong}
	For all $\epsilon > 0$, there is a constant $C(\epsilon) > 0$
	such that
	\[
		\verts{\Delta_E} \leq C(\epsilon) N_E^{6 + \epsilon},
	\]
	for any elliptic curve $E$ over $\Q$ with minimal discriminant
	$\Delta_E$ and conductor $N_E$.
\end{conjecture}

\begin{conjecture}[{Modified Szpiro's conjecture \cite[Conjecture~4]{oesterle}}]
	\label{conj:modified-szpiro}
	For each $\epsilon > 0$, there is a constant $C(\epsilon) > 0$ such
	that
	\[
		\max\set{\verts{c_4}^3, \verts{c_6}^2} \leq C(\epsilon) N_E^{6 + \epsilon},
	\]
	for any elliptic curve $E$ over $\Q$ with invariants $c_4$,
	$c_6$, and conductor $N_E$.
\end{conjecture}

The many arithmetic implications of Szpiro's conjecture
(and its explicit variations)
famously include
the Mordell conjecture \cite{elkies},
Fermat's Last Theorem \cite{granville-tucker},
Baker's theorem \cite{granville-tucker},
Roth's theorem on Diophantine approximations \cite{granville-tucker},
Lang's height conjecture \cite[p.~420]{hindry-silverman-1988},
and the non-existence of Siegel zeroes for certain
L-functions \cite[Theorem~2]{granville-stark}.

\begin{remark}
	It is worth noting the significant attention that
	Szpiro's conjecture has received since the 2012 release of Shinichi
	Mochizuki's four preprints claiming a proof via
	inter-universal Teichm\"{u}ller theory.
	These four papers have since been published
	\cite{mochizuki1, mochizuki2, mochizuki3, mochizuki4}
	with an additional follow-up article \cite{mochizuki-explicit},
	although the academic disagreements
	have not yet been completely resolved to the author's knowledge
	(c.f. \cite{scholze-stix, scholze-quanta, scholze-nature, dupuy-hilado, mochizuki-essential}).
\end{remark}

\subsection{\texorpdfstring{$abc$}{abc} conjecture}

Seeking to formulate Szpiro's conjecture without referring
to elliptic curves, Masser and
Oesterl\'{e} formulated the $abc$ conjecture in 1985
(cf. \cite[Section~3]{oesterle}).

\begin{conjecture}[{The $abc$ conjecture \cite[Conjecture~3]{oesterle}}]
	\label{conj:abc}
	For every positive real number $\epsilon$, there exists a positive
	real number $C(\epsilon)$ such that
	\[
		c < C(\epsilon) \rad(abc)^{1 + \epsilon},
	\]
	for every triple $(a, b, c)$ of coprime positive integers such that
	$a + b = c$.
\end{conjecture}

The $abc$ conjecture was then shown to be equivalent to the
modified Szpiro's conjecture by Oesterl\'{e} and Nitaj
(cf. \cite[Section~3]{oesterle}, \cite[Section~4]{goldfeld1999},
\cite[p.~1769]{notices-szpiro}).
One way to see that Szpiro's conjecture implies a
weak $abc$ conjecture is to take the Hellegouarch--Frey curve
\[
	E_{a, b}: y^2 = x(x - a)(x - b),
\]
for any positive coprime integers $(a, b, c)$ such that
$a + b = c$.
The elliptic curve $E_{a, b}$, after taking
a minimal model,
has minimal discriminant $\Delta_{a, b} = 2^{-s} (abc)^2$ and
conductor $N_{a, b} = 2^{-t} \rad(abc)$ where $s$ and $t$
are bounded integers. If Conjecture~\ref{conj:modified-szpiro}
holds for $E_{a, b}$, then for all $\epsilon > 0$,
there is a constant $C(\epsilon) > 0$ such that
\begin{align*}
	\verts{\Delta_{a, b}} &\leq C(\epsilon) N_{a, b}^{6 + \epsilon}, \\
	\verts{2^{-s} (abc)^2} &\leq C(\epsilon) \paren{2^{-t} \rad(abc)}^{6 + \epsilon}, \\
	(abc)^2 &\leq C(\epsilon) \frac{\paren{2^{-t}}^{6 + \epsilon}}{\verts{2^{-s}}}
		\rad(abc)^{6 + \epsilon}.
\end{align*}
Folding $\frac{\paren{2^{-t}}^{6 + \epsilon}}{\verts{2^{-s}}}$ into
$C(\epsilon)$ and observing that $c = a + b$ is greater than $a$ and
$b$, we obtain
\begin{align*}
	c^4 \leq (abc)^2 &\leq C(\epsilon) \paren{\rad(abc)}^{6 + \epsilon} \leq C(\epsilon) \rad(abc)^{6 + 6\epsilon},
\end{align*}
and thus,
\begin{align*}
	c &\leq C(\epsilon) \rad(abc)^{\frac{3}{2} + \epsilon}.
\end{align*}

\begin{remark}
	In fact, Baker~\cite{baker} and
	Laishram--Shorey~\cite{laishram-shorey} give empirical evidence
	for an explicit version of Conjecture~\ref{conj:abc}
	with $\epsilon = \frac{3}{4}$ and $C\paren{\frac{3}{4}} \leq 1$
	(these exponents can be even further optimized,
	cf. \cite{chim-shorey-sinha}).
\end{remark}

\section{\texorpdfstring{$abcd\ldots$}{abcd...} conjectures}
\subsection{Background and notation}
There are generalizations of the $abc$ conjecture in various
directions, such as the uniform version over general number fields
given by Vojta~\cite{vojta-1987}. Vojta~\cite{vojta-1998} also
proposed a higher-dimensional generalization over general number
fields using the arithmetic truncated counting function in the
language of Arakelov theory. We try to present the language of
\cite{vojta-1998} (e.g.~using the arithmetic model) as well as the
language of \cite{looper1, looper2} (e.g.~using
local height functions).

Let $K$ be a number field or one-variable function field of
characteristic $0$ with set of places $M_K$.
If $K$ is a function field, then $K = k(\M_K)$ (the function field
of $\M_K$) for some smooth
curve $\M_K$ over a field of constants $k_0 \subset K$ and
there is a bijection between $M_K$ and the closed points of
$\M_K$.
If $K$ is a number field, then there is an arithmetic scheme
$\M_K$ whose closed points are also in bijection with $M_K$
and whose non-archimedean part $(\M_K)_\na$ is just $\spec \sco_K$.
An arithmetic model $\scx$ of $X$ is an integral, flat, separated
scheme $\scx_\na$ of finite type over $(\M_K)_\na$
(plus additional archimedean information) with an isomorphism
$X \cong \scx_\na \times_{(\M_K)_\na} K$.
\begin{remark}
	In this article, we will largely ignore the Green functions
	and archimedean phenomena (cf. \cite{vojta-1998}).	
\end{remark}

For a place $\nu \in M_K$, denote the residue field of $K$ at $\nu$
by $k_\nu$.
Define $\mu_\nu$ and $r_\nu$ via
\begin{align*}
	\mu_\nu :&= \begin{cases}
			\log\paren{\# k_\nu}
				&\textrm{ if } K \textrm{ is a number field,}\\
			[k_\nu:k_0]
				&\textrm{ if } K \textrm{ is a function field.}
		\end{cases} \\
	r_\nu :&= \begin{cases}
			\frac{[K_\nu:\Q_\nu]}{[K:\Q]}
				&\textrm{ if } K \textrm{ is a number field,}\\
			\frac{[K_\nu:k_0(t)_\nu]}{[K:k_0(t)]}
				&\textrm{ if } K \textrm{ is a function field.}
		\end{cases}
\end{align*}
For a Cartier divisor $D$ on $\M_K$,
let $n_\nu(D)$ be the multiplicity of $D$ at the closed point of
$\M_K$ corresponding to $\nu$.

Relative to a divisor $D$ of a smooth complete variety $X$
over a field $K$ and $\nu$-adic local heights
(i.e.~local Weil functions) $\lambda_{D, \nu}$
(cf. \cite[Section~1.3]{vojta-1987},
\cite[Section~B.8]{hindry-silverman-book}),
we will denote the height of a point $P \in X(\overline{K})
\backslash \supp(D)$ by
\[
	h_D(P) := \sum_{\nu \in M_K} r_\nu \lambda_{D, \nu}(P).
\]
Alternatively, the height can be given as
\[
	h_D(P) = \frac{\deg \sigma^* D}{[L:K]}
\]
where $D$ is a Cartier divisor on $\scx$,
$P$ is not in the support of $D$,
$E$ is a finite field extension of $K$ containing
$K(P)$, and $\sigma$ is the map $\M_L \rightarrow \scx$
corresponding to $P$ (cf. \cite[Section~1]{vojta-1998}.

Following Vojta~\cite[Section~1]{vojta-1998},
we define the truncated counting function from
higher-dimensional Nevanlinna theory (cf.
\cite[Section~3.4]{vojta-1987}, \cite[Section~2.2]{looper1}).
\begin{definition}
	Let $X$ be a smooth complete variety over a number field 
	or one-variable function field $K$, and
	$S$ be a finite set of places of $K$ containing
	the archimedean places.
	For an effective Cartier divisor $D$ on $X$ and a point
	$P \in X(\overline{K}) \backslash \supp D$, define the
	\emph{truncated counting function} by
	\[
		N^{(1)}_{K, S} (D, P) := \frac{1}{[K(P):K]}
			\sum_{\substack{\nu \in M_{K(P)} \\ \nu \not \mid S}}
			\min\set{1, n_\nu(D)} \mu_\nu
	\]
\end{definition}
Relative to the local height functions, the truncated counting
function can be given as
\[
	N^{(1)}_{K, S} (D, P) := \frac{1}{[K(P):K]}
			\sum_{\substack{\nu \in M_{K(P)} \\ \nu \not \mid S}}
			\chi_{\R_{>0}}\paren{\lambda_{D, \nu}(P)} \mu_\nu,
\]
where $\chi_{\R_{>0}}$ is the characteristic function
for the positive real numbers.

We also recall the definition of a normal crossings divisor (cf.
\cite[Section~1]{vojta-1998},
\cite[Definition~4.1.1.]{lazarsfeld},
\cite[\href{https://stacks.math.columbia.edu/tag/0CBN}{Tag 0CBN}]{stacks-project}).
\begin{definition}[Normal crossings divisor]
	A Cartier divisor $D$ on a smooth variety $X$
	is a \emph{normal crossings divisor} if $D$ is a formal sum $\sum_i D_i$ of
	distinct irreducible subvarieties $D_i$ and if it can be represented
	for every point $P \in X$ by a principal divisor
	$(x_1 \cdots x_r)$ in the completed local ring
	$\hat{\sco}_{P, X}$, where $x_1, \ldots, x_r$ is part of a regular
	sequence for $\hat{\sco}_{P, X}$.
\end{definition}

Finally, we define the logarithmic discriminant of a finite field
extension.
Here we follow Looper~\cite[Section~2.2]{looper1},
but a uniform definition
can also be given by observing that there is a finite morphism
\[
	(\M_L)_\na \rightarrow (\M_K)_\na,
\]
and then taking the degree of its ramification divisor
(viz. \cite[Section~1]{vojta-1998}).
\begin{definition}[Logarithmic discriminant]
	Let $L/K$ be a finite extension of a number field or one-variable
	function field of characteristic zero,
	let $D_{L/K}$ be its discriminant, and
	if $L/K$ is a function field then let $g(L)$ denote its genus.
	Define the \emph{logarithmic discriminant} of $L/K$ to be
	\[
		d_{L/K} := \begin{cases}
			\frac{\log \paren{\sum_{\nu \in (\M_K)_\na} \verts{D_{L/K}}_\nu^{-r_\nu}}}{[L:K]}
				&\textrm{ if } L/K \textrm{ is a number field,}\\
			\frac{2g(L) - 2}{[L:K]}
				&\textrm{ if } L/K \textrm{ is a function field.}
		\end{cases}
	\]
\end{definition}

\begin{remark}
	Looper~\cite{looper1} refers to
	Kim--Thakur--Voloch~\cite[Section~2]{kim-thakur-voloch} for the
	definition of the genus of a function field in one variable.
	It can be defined instead in terms of the genus of $\M_K$
	or in terms of divisors and the Riemann--Roch theorem
	(viz. \cite[Section~2.1]{chevalley},
	\cite[Section~I.2]{lang-functions}, or
	\cite[Definition~1.4.15]{stichtenoth}).
\end{remark}

\subsection{Vojta's \texorpdfstring{$abcde\ldots$}{abcde...} conjecture}
With the $abc$ conjecture, a natural question is to ask whether
one might expect a similar inequality between the maximum of an
$n$-tuple of positive coprime integers and their radical.
Using a higher-dimensional version of Nevanlinna theory,
Vojta~\cite{vojta-1998} posed a higher-dimensional generalization
of the $abc$ conjecture called the ``$abcde\ldots$ conjecture''.
\begin{conjecture}[{The $abcde\ldots$ conjecture \cite[Conjecture~2.3]{vojta-1998}}]
	\label{conj:abcde}
	Let $X$ be a smooth complete variety over a number field or
	one-variable function field $K$ of
	characteristic $0$, $S$ be a finite set of places of $K$ containing
	the archimedean places, $D$ be a normal crossings divisor on $X$,
	$K_X$ be the canonical divisor on $X$, $\sca$ be a big line
	bundle on $X$, $r$ be a positive integer, and $\epsilon > 0$.
	Then there exists a proper Zariski-closed subset
	$\calz = \calz(K, S, X, D, \sca, r, \epsilon) \subsetneq X$ such that
	\[
		N^{(1)}(D, P) \geq h_{K_X + D}(P) - \epsilon h_\sca(P) - d_{K(P)/K} + O(1),
	\]
	for all $P \in X(\overline{K}) \backslash
	\calz$ with $[K(P):K] \leq r$.
\end{conjecture}

Conjecture~\ref{conj:abcde} in this original form is known to
be false, but can be salvaged with modification as we
will mention later.
First, we describe how it generalizes the $abc$ conjecture
(Conjecture~\ref{conj:abc}).
A triple $(a, b, c)$ of coprime positive integers such that $a + b = c$
corresponds to a point $P = [a:b:-c]$ on the line
$X := \set{Z_0 + Z_1 + Z_2 = 0} \cong \P_\Q^1$ in
$\P_\Q^2$. Let $\sca := \sco(1)$ and let $D$
be the normal crossings divisor
$\set{Z_0 = 0} + \set{Z_1 = 0} + \set{Z_2 = 0}$
(in the coordinates of the ambient $\P_\Q^2$).
Then $h_\sca(P) = \log(c)$ and $N^{(1)}(D, P) = \rad(abc)$.
Notice that $h_{K_X + D}(P)
= \log(c) + O(1)$ since $\omega_X \cong \sco(-2)$,
$D$ corresponds to three distinct points
on $X$, and $\sco(K_X + D) \cong \sco(1)$.
In this example with $S = \set{\infty}$ and $r=1$,
the inequality of Conjecture~\ref{conj:abcde} becomes
\begin{align*}
	\rad(abc) &\geq \log(c) + O(1) - \epsilon \log(c) - 0 + O(1) \\
		&\geq (1 - \epsilon) \log(c) + O(1).
\end{align*}
With a different $\epsilon$
(such as $\epsilon \mapsto 1 - \frac{1}{1+\epsilon}$),
\[
	\log(c) \leq (1 + \epsilon) \rad(c) + O(1).
\]
Exponentiating yields the desired inequality of the $abc$ conjecture,
with a constant $C$ not necessarily depending on $\epsilon$.

Repeating the same procedure any $n$ yields an analogous conclusion
for generic $n$-tuples
$(z_0, \ldots, z_{n-1})$ of coprime positive integers such
that $\sum_{i=0}^{n-2} z_i = z_{n-1}$.
The $abcde\ldots$ conjecture implies that, outside of a proper
Zariski-closed subset of the hyperplane
$X := \set{\sum_{i=0}^{n-1} Z_i = 0} \cong \P_\Q^{n-1}$,
all such tuples satisfy the inequality
\begin{align}
	\label{eq:n-abcde}
	z_{n-1} \leq C(\epsilon) \rad\paren{\prod_{i=0}^{n-2} z_i}^{1 + \epsilon}.
\end{align}
The nomenclature of the ``$abcde\ldots$ conjecture'' for
Conjecture~\ref{conj:abcde} is evident in
this consequence, if one were to relabel $(z_0, \ldots, z_{n-1})$
as $(a, b, c, d, e, \ldots)$.

For $n=4$, inequality~\ref{eq:n-abcde} is the $abcd$ conjecture
described by Granville~\cite[Section~8]{granville}
(cf. \cite[Chapter~IV.3]{princeton}).
\begin{conjecture}[{The $abcd$ conjecture, version $1$ \cite[Section~8]{granville}}]
	\label{conj:abcd-granville}
	For every positive real number $\epsilon$, there exists a positive
	real number $C(\epsilon)$ such that
	\[
		d < C(\epsilon) \rad(abcd)^{1 + \epsilon},
	\]
	for every quadruple $(a, b, c, d)$ of coprime positive integers for
	which $a + b + c = d$ outside of a Zariski-closed subset.
\end{conjecture}
\begin{remark}
	Granville~\cite[Section~8]{granville} describes another variant of
	Conjecture~\ref{conj:abcd-granville},
	removing the ``outside of a Zariski-closed subset'' condition with
	a larger exponent in the inequality:
	\[
		d < C(\epsilon) \rad(abcd)^{3 + \epsilon}.
	\]
	This variation of the Conjecture~\ref{conj:abcd-granville}
	is the case $n=4$ of the $n$-conjecture of
	Browkin--Brzezinski~\cite[Section~2]{browkin-brzezinski}
	(also sometimes called the $abcd$ conjecture,
	e.g. \cite[Section~4]{ramaekers}, \cite[Section~2]{voloch}),
	which asserts a general exponent of $2n - 5 + \epsilon$.
\end{remark}

Another version of the $abcd$ conjecture is described by
Looper~\cite[Conjecture~2.1]{looper1}, using
a more general form of inequality~\ref{eq:n-abcde}
in the language of heights. Define
\[
	\rad(P) = \rad\paren{\sbrac{z_0:\ldots:z_{n-1}}} :=
		\frac{1}{[K(P):K]}
		\sum_{\substack{\nu \in (M_K)_\na \\ \nu(z_i) \neq \nu(z_j)
		\textrm{ for some } i, j}} \mu_\nu.
\]

\begin{conjecture}[{The $abcd$ conjecture, version $2$ \cite[Conjecture~2.1]{looper1}}]
	\label{conj:abcd}
	Let $K$ be a number field or a $1$-dimensional function field of
	characteristic $0$. Let $n \geq 3$, $[Z_1, \ldots, Z_n]$ be the
	standard coordinates on $\P^{n-1}_K$, and $\calh \subset \P^{n-1}_K$
	be the hyperplane given by $\sum_i Z_i = 0$. For any $\epsilon > 0$,
	there is a proper Zariski-closed subset
	$\calz = \calz(K, \epsilon, n) \subsetneq \calh$ and a constant
	$C_{K, \calz, \epsilon, n}$ such that
	\[
		h(P) < (1 + \epsilon) \rad(P) + C_{K, \calz, \epsilon, n},
	\]
	for all $P \in \calh \backslash \calz$.
\end{conjecture}

\begin{remark}
	For the remainder of this article, ``$abcd$ conjecture''
	will mean the Version $2$ (Conjecture~\ref{conj:abcd})
	as described by Looper~\cite{looper1}.
\end{remark}

The original form of Vojta's $abcde\ldots$ conjecture,
given in Conjecture~\ref{conj:abcde}, was shown
to be false by Masser~\cite{masser}.
Looper~\cite[Section~2.2]{looper1}
notes that it can be replaced by the following weaker conjecture
(with $P \in X(K)$ and without the logarithmic discriminant).

\begin{conjecture}[{The weak $abcde\ldots$ conjecture \cite[Conjecture~2.2]{looper2}}]
	\label{conj:abcde-weak}
	Let $X$ be a smooth complete variety over a number field or
	one-variable function field $K$ of
	characteristic $0$, $S$ be a finite set of places of $K$ containing
	the archimedean places, $D$ be a normal crossings divisor on $X$,
	$K_X$ be the canonical divisor on $X$, $\sca$ be a big line
	bundle on $X$, and $\epsilon > 0$.
	Then there exists a proper Zariski-closed subset
	$\calz = \calz(K, S, X, D, \sca, r, \epsilon) \subsetneq X$ such that
	\[
		N^{(1)}(D, P) \geq h_{K_X + D}(P) - \epsilon h_\sca(P) + O(1),
	\]
	for all $P \in X(K) \backslash \calz$.
\end{conjecture}

This form is only weaker by ignoring the logarithmic
discriminant and requiring that $P$ is $K$-rational rather than be of
bounded degree over $K$. It still implies the $abc$ conjecture and the
$abcd$ conjecture.

\section{Uniform boundedness}

\subsection{Uniform boundedness of torsion points}
In the first decade of the 20th century, Levi conjectured
a classification of torsion groups for elliptic curves
over $\Q$ (cf. \cite[Section~6]{schappacher-schoof})
that was later eventually proved by Mazur~\cite{mazur1, mazur2}
and then extended to general number fields by Kamienny~\cite{kamienny},
Kamienny--Mazur~\cite{kamienny-mazur}, and Merel~\cite{merel}.
\begin{theorem}[{Torsion theorem \cite{merel}}]
	\label{thm:torsion}
	Fix a positive integer $D$. There is a positive integer
	$C(D)$ such that for all elliptic curves $E$ defined over a
	number field $K$ of degree $D$,
	the number of $K$-rational torsion points on $E$ is uniformly
	bounded by $C(D)$:
	\[
		\#E(K)_{\mathrm{tors}} \leq C(D).
	\]
\end{theorem}

\begin{remark}
	In fact, the uniform bound of Theorem~\ref{thm:torsion}
	can be made effective by the results of Parent~\cite{parent}.
\end{remark}

The natural generalization of the torsion theorem to abelian varieties
is an open problem for dimension greater than $1$.
\begin{conjecture}[{Uniform boundedness conjecture for torsion points of abelian varieties}]
	\label{conj:torsion}
	Fix positive integers $g$ and $D$. There is a positive integer
	$C(g, D)$ such that for all abelian varieties $A$ of dimension $g$
	defined over a number field $K$ of degree $D$,
	the number of $K$-rational torsion points on $A$ is uniformly
	bounded by $C(g, D)$:
	\[
		\#A(K)_{\mathrm{tors}} \leq C(g, D).
	\]
\end{conjecture}

\noindent
The general uniform boundedness conjecture for torsion points of
abelian varieties is wide open. However,
it is sufficient to prove the conjecture for Jacobian
varieties due to Cadoret--Tamagawa~\cite{cadoret-tamagawa-note};
they show that Conjecture~\ref{conj:torsion} for abelian
varieties of dimension $g$ follows from the statement of
Conjecture~\ref{conj:torsion} for Jacobian varieties
of dimension $1 + 6^{8g}(8g-1)!\frac{8g(8g-1)}{2}$.
They also demonstrated a version of uniform boundedness for the
$p$-primary part of torsion in families of $g$-dimensional
abelian varieties parametrized by curves
\cite[Theorem~1.1]{cadoret-tamagawa-ubc},
\cite[Corollary~4.3]{cadoret-tamagawa-ubc2}.

Otherwise, there has been partial progress in specific cases,
such as for certain abelian surfaces of CM-type
(cf. \cite[Corollaries~2-3]{silverberg},
\cite[Th\'{e}or\`{e}me~2]{gaudron-remond}) and
for abelian varieties with specific anisotropic
reduction (cf. \cite[Corollary~2 and Theorem~4]{clark-xarles}).
In a parallel to the development of the $abcd$ conjecture,
it is known due to Clark--Xarles~\cite[Section~6]{clark-xarles}
that a higher-dimensional analogue of the
weak Szpiro's conjecture
implies the uniform boundedness of torsion points
for Hilbert--Blumenthal abelian varieties.

\subsection{Uniform boundedness of preperiodic points}
In the simplest dynamical setting, we replace endomorphisms of
elliptic curves with endomorphisms of projective space $\P^N$.
\begin{theorem}[{Northcott's theorem \cite{northcott}}]
	Fix positive integers $N \geq 1$ and $d \geq 2$. For all
	degree $d$ morphisms $f: \P^N \rightarrow \P^N$
	defined over a number field $K$, the number of
	$K$-rational preperiodic points is finite.
\end{theorem}

A uniform boundedness conjecture for preperiodic points
of an endomorphism on projective space was posed by 
Morton--Silverman~\cite{morton-silverman}.
\begin{conjecture}[{Dynamical uniform boundedness conjecture \cite{morton-silverman}}]
	\label{conj:dynamical-uniform}
	Fix integers $N \geq 1$, $d \geq 2$, and $D \geq 1$. There is a
	positive integer $C(N, d, D)$ such that for all degree $d$ morphisms
	$f: \P^N \rightarrow \P^N$ defined over a number field $K$
	of degree $D$, the number of $K$-rational
	preperiodic points is uniformly bounded by $C(N, d, D)$:
	\[
		\# \preper(f, \P^N_K) \leq C(N, d, D).
	\]
\end{conjecture}

Fakhruddin~\cite[Corollary~2.4]{fakhruddin} showed that
the dynamical uniform boundedness conjecture implies the
uniform boundedness conjecture for torsion points on
abelian varieties.

Assuming the $abcd$ conjecture, Looper
proved a version of the dynamical uniform boundedness conjecture
of Morton--Silverman~\cite{morton-silverman},
first for unicritical polynomials \cite{looper1}
and then for general single-variable
polynomials \cite{looper2}.
\begin{theorem}[{\cite[Theorem~1.2]{looper1, looper2}}]
	\label{thm:looper-uniform}
	Fix an integer $d \geq 2$.
	Let $K$ be a number field or a $1$-dimensional function field
	of characteristic $0$.
	If the $abcd$ conjecture holds, then there is a positive integer
	$C(d, K)$ such that for any polynomial $f \in K[z]$ (that is
	furthermore isotrivial if $K$ is a function field) of degree $d$,
	the number of $K$-rational preperiodic points of $f$ is
	uniformly bounded by $C(d, K)$:
	\[
		\# \preper(f, \P^1_K) \leq C(d, K).
	\]
\end{theorem}
The proof of Theorem~\ref{thm:looper-uniform} uses several
ingredients. A key ingredient is a global equidistribution statement
pieced together from local $\nu$-adic Julia-Fatou equidistribution
results obtained from the potential theory
of algebraically closed and complete non-Archimedean metrized fields
via the Berkovich projective line
(for more on this general area, see \cite{baker-rumely}).
Besides the input
of height machinery and bounds on heights of preperiodic points
of polynomial functions
(e.g. from \cite{ingram-lower, looper-lower}),
one of the other main ideas is to show that the prime factors of
differences of preperiodic points
$z_i - z_j$ are typically in the places of bad reduction when
$f$ has many preperiodic points. The $abcd$ conjecture
(Conjecture~\ref{conj:abcd}) can then
be applied to combinatorial ``polygons'' constructed from
preperiodic points giving points on the projective
hyperplanes $\calh$.

Following Looper~\cite{looper1} with some case-by-case calculations,
Panraksa~\cite{panraksa} showed that the usual
$abc$ conjecture implies that there are no rational non-fixed
periodic points for large-degree unicritical polynomials.
\begin{theorem}[{\cite[Theorem~3]{panraksa}}]
	\label{thm:panraksa-unicritical}
	If the $abc$ conjecture holds, then
	the number of rational preperiodic points of $f$ is
	bounded,
	\[
		\# \preper(f_{d, c}, \P^1_\Q) \leq 4,
	\]
	for all unicritical polynomials
	$f_{d, c}(z) := z^d + c \in \Q[z]$ with $d$ sufficiently large.
	Furthermore, if $c \neq -1$ then
	$f_{d, c}$ has no rational periodic points of exact period greater
	than $1$ for $d$ sufficiently large.
\end{theorem}

The idea of the proof of Theorem~\ref{thm:panraksa-unicritical}
follows the difference-of-preperiodic-points idea used
in \cite{looper1, looper2}. The main observation is that
by the $abc$ conjecture for the positive integer triple
$\paren{\verts{Z_1^d}, \verts{Z_2^d}, \verts{(Z_3 - Z_2) Z^{d-1}}}$,
the particular system of equations
\begin{align*}
	Z_2^d - Z_1^d &= (Z_3 - Z_2) Z^{d-1} \neq 0, \\
	\gcd(Z_1, Z_2) &= 1, \\
	\max \set{\verts{Z_1}, \verts{Z_2}, \verts{Z_3}} &= Z_3, \\
	\max \set{\verts{Z_1}, \verts{Z_2}, \verts{Z_3}, \verts{Z}} &> 1,
\end{align*}
has no integral solutions $(Z_1, Z_2, Z_3, Z)$
for sufficiently large $d$. Periodic points of $f_{d, c}$
generate such systems of equations by expressing elements of
an orbit as
\[
	\set{z_1 = \frac{Z_1}{Z}, \ldots, z_n = \frac{Z_N}{Z}},
\]
and looking at the differences, e.g.~$z_3 - z_2 = z_2^d - z_1^3$.

The statement of Theorem~\ref{thm:panraksa-unicritical} 
about the non-existence of rational periodic points
for large enough $d$ can be
partially lifted from $\Q$ to quadratic number fields $K$
via the following observation. For a morphism
$f: \P^1 \rightarrow \P^1$, $f^{(j)} := f \circ f^{(j-1)}$ is
the $j$-th iterate of $f$.
\begin{lemma}[{\cite[Lemma~3.3]{zhang}}]
	\label{lem:rational-points}
	Let $K$ be a quadratic number field
	and $f_{d, c}(z) := z^d + c \in \Q[z]$ be a unicritical
	polynomial of degree $d > 1$.
	Let $\set{z_0, \ldots z_{N-1}} \subset K$ be
	an orbit of periodic points of $f_{d, c}$ of odd period $N > 1$.
	If there is a positive integer
	$j_i < N$ and a nontrivial $\sigma_i \in \gal(K/\Q)$ for each
	$z_i$ such that
	\[
		\sigma_i (z_i) = f_{d, c}^{(j_i)}(z_i),
	\]
	then each $z_i$ is a rational number.
\end{lemma}

\begin{remark}
	The condition of Lemma~\ref{lem:rational-points} is
	called the ``Galois--dynamics correspondence''
	\cite[Definition~1.1]{zhang}. This condition
	is satisfied, for example, when the dynatomic
	polynomial $\Phi_N(z) \in \Q[z]$ of $f_{d, c}$
	is irreducible \cite[Proposition~1.5]{zhang}.
\end{remark}

\begin{corollary}
	\label{cor:panraksa-quad}
	Let $K$ be a quadratic number field and $N$ be an odd integer
	greater than $1$.	
	If the $abcd$ conjecture holds, the Galois--dynamics
	correspondence holds for every period-$N$ orbit of $f_{d, c}$ in
	$K$, and $c \neq -1$, then
	$f_{d, c}$ has no $K$-rational periodic points of exact period
	greater than $1$ for $d$ sufficiently large.
\end{corollary}

Using Theorem~\ref{thm:panraksa-unicritical},
Panraksa observed that the uniform
bound of Theorem~\ref{thm:panraksa-unicritical} for $K=\Q$ is actually
an absolute constant not dependent on $d$.
\begin{corollary}[{\cite[Theorem~4]{panraksa}}]
	\label{cor:panraksa-absolute}
	If the $abcd$ conjecture holds, then there is a positive integer
	$C$ such that for all unicritical polynomials
	$f_{d, c}(z) := z^d + c \in \Q[z]$ with $d \geq 2$,
	\[
		\# \preper(f_{d, c}, \P^1_\Q) \leq C.
	\]
\end{corollary}

While the uniform bound on the number of preperiodic points of
unicritical polynomials in Corollary~\ref{cor:panraksa-absolute}
is for the rational numbers,
there is an unconditional bound on the possible periods
of periodic points for number fields $K$ due
to a result of Morton--Silverman~\cite[Corollary~B]{morton-silverman}.
Their result for periodic points of endomorphisms over number fields
was then extended to preperiodic points over all global fields by
Canci--Paladino~\cite{canci-paladino}. We state the number field
version of Canci--Paladino's result
for the forward orbit $\sco_f(P) = \set{f^{(j)}(P) \mid j \in \N}$
of a preperiodic point $P$ of an endomorphism $f$.
\begin{theorem}[{\cite[Theorem~1]{canci-paladino}}]
	\label{thm:morton-silverman}
	Let $K$ be a number field of degree $D$, $S$ be a finite set
	of places of $K$, and $f$ be an endomorphism of $\P^1_K$
	of degree $d$ defined over $K$ with good reduction outside
	of $S$. If $P \in \P^1_K$ is a preperiodic point of $f$
	and $N := \verts{\sco_f(P)}$, then
	\[
		N \leq \max \set{
			\paren{2^{16\verts{S} - 8} + 3} \paren{12\verts{S} \log\paren{5\verts{S}}}^D,
			\paren{\paren{12\verts{S} + 24} \log\paren{5\verts{S}+5}}^{4D}}.
	\]
\end{theorem}

Specializing to periodic points of quadratic rational maps $f$,
Canci~\cite{canci} used Morton--Silverman's bound to show
that quadratic rational maps
with good reduction outside a set of places $S$ typically
do not have periodic points of large period.
\begin{theorem}[{\cite[Theorem~1']{canci}}]
	Let $K$ be a number field with ring of integers $\sco_K$
	and let $S$ be a finite set of places of $K$.
	There are only finitely many $\pgl_2\paren{\sco_{K, S}}$-conjugacy
	classes of quadratic rational maps defined over $K$ with good
	reduction outside of $S$ and with a periodic point of minimal period
	greater than $3$.
\end{theorem}

For periodic points of unicritical polynomials with coefficients in $\sco_K$,
Panraksa observed in an unpublished note \cite{panraksa-unpublished}
that Theorem~\ref{thm:morton-silverman}
implies a uniform bound on possible periods.
\begin{corollary}
	\label{cor:panraksa-morton-silverman}
	Let $K$ be a number field of degree $D$
	with ring of integers $\sco_K$.
	For the unicritical polynomial
	$f_{d, c}(z) := x^d + c \in \sco_K[z]$ with $d > 1$,
	there is a positive integer $C_D$ depending only
	on the degree $D$ such that $f_{d, c}$
	has no periodic points of period $N \geq C_D$ in $K$.
\end{corollary}
\begin{proof}
	Take $S$ to be the set of archimedean places in $K$.
	Since $c \in \sco_K$, the unicritical polynomial $f_{d, c}$
	has good reduction everywhere outside of $S$.
	Furthermore, $\verts{S} = r_K + s_K \leq [K:\Q]$.
	By Theorem~\ref{thm:morton-silverman},
	$N \leq C_D$ where
	\[
		C_D = \paren{12\paren{D + 1}
			\log \paren{5\paren{D+1}}}^{4[K:\Q]}.
	\]
\end{proof}

\subsection{Uniform boundedness of small points}
One of the well-known consequences of Szpiro's conjecture
due to Hindry--Silverman~\cite[p.~420]{hindry-silverman-1988}
is Lang's height conjecture about a lower bound for the
N\'{e}ron--Tate (canonical) height.
\begin{conjecture}[{Lang's height conjecture \cite[p.92]{lang-elliptic}}]
	\label{conj:lang-height}
	Let $K$ be a number field.
	There is a positive constant $C(K)$ such that
	for any elliptic curve $E$ defined over $K$ with minimal
	discriminant $\Delta_E$, the
	N\'{e}ron--Tate height of every non-torsion point $P \in E(K)$
	satisfies
	\[
		\hat{h}(P) \geq C(K) \log \paren{N_{K/\Q} \Delta_E}.
	\]
\end{conjecture}
\begin{remark}
	If $E/K$ has good reduction everywhere, then $N_{K/\Q} \Delta_E = 1$
	and the bound in Lang's height conjecture is trivial.
	In this case, Hindry--Silverman~\cite[Corollaire~2]{hindry-silverman-1999}
	obtained a nontrivial lower bound for the N\'{e}ron--Tate height
	of non-torsion points $P \in E(K)$ solely in terms of $D := [K:\Q]$:
	\[
		\hat{h}(P) \geq \paren{10^{18} D^3 (\log(D))^2 }^{-1}.
	\]
\end{remark}

For an endomorphism $\phi$ of $\P^N$ of degree $d \geq 2$,
the canonical height
$\hat{h}_\phi$ is the unique real-valued function on
$\P^N(\overline{\Q})$ such that $\hat{h}_\phi(P) = h(P) + O(1)$
and $\hat{h}_\phi(\phi(P)) = d\hat{h}_\phi(P)$,
where $h(P)$ is the Weil height. It is similarly characterized by
the property that the height $\hat{h}_\phi(P)$ vanishes if and only
if $P$ is a preperiodic point for
$\phi$ \cite[Theorem~3.22]{silverman-dynamics}.
Silverman~\cite[Conjecture~4.98]{silverman-dynamics} formulated
a dynamical version of Lang's conjecture (not to be confused with
the dynamical Mordell--Lang conjecture) that captures how
the canonical height $\hat{h}_f$ describes the ``non-preperiodicity''
of a point $P$ relative to a dynamical system.

Let $\rat_d$ be the moduli space of rational maps of degree $d$.
The quotient variety $\scm_d := \rat_d / \pgl_2$ by the $\pgl_2$
conjugacy action is a moduli space classifying degree $d$
dynamical systems on $\P^1$
(cf. \cite[Section~4.4]{silverman-dynamics}).
Fix a projective embedding $\scm_d \hookrightarrow \P^n$ so
that a rational function $\phi \in K(z)$ has a corresponding
Weil height from $\abrac{\phi} \in \scm_d(K)$. Let
the minimal resultant of $\phi$ be
\[
	\mfR_\phi := \prod_\mfp \mfp^{\epsilon_\mfp(\phi)},
\]
where $\epsilon_\mfp(\phi)$ is the greatest exponent of $\mfp$
dividing the resultant of the $\pgl_2$-conjugate of $\phi$ with
the best reduction at $\mfp$
(cf. \cite[Section~4.11]{silverman-dynamics}). This allows for
the accounting of twists of $\phi$ which have the same $\scm_d$
class but do not necessarily have the same Weil height.

\begin{conjecture}[{Dynamical Lang conjecture \cite[Conjecture~4.98]{silverman-dynamics}}]
	Fix an embedding of the moduli space $\scm_d$ in projective space
	and let $h_{\scm_d}$ denote the associated height function.
	Let $K$ be a number field and $d \geq 2$ be an integer.
	Then there is a positive constant $C(K, d)$ such that
	\[
		\hat{h}_\phi(P) \geq C(K, d)
			\max\set{\log (N_{K/\Q} \mfR_\phi), h_{\scm_d}\paren{\abrac{\phi}}}.
	\]
	for all rational maps $\phi \in K(z)$ of degree $d$ and all
	non-preperiodic points $P \in \P^1_K$.
\end{conjecture}

Using the $abcd$ conjecture along the same lines as
the proof of Theorem~\ref{thm:looper-uniform},
Looper gives a weaker version of the
dynamical Lang conjecture.
Here, this version uses the critical height, which
is the sum of canonical heights at critical points
\[
	\hat{h}_\crit (\phi) := \sum_{P \in \Crit(\phi)} \hat{h}_\phi(P),
\]
and is commensurate with the Weil height $h_{\scm_d}$
\cite[Theorem~1]{ingram}.

\begin{theorem}[{\cite[Theorem~1.3]{looper2}}]
	\label{thm:looper-lang}
	Fix an integer $d \geq 2$.
	Let $K$ be a number field or a $1$-dimensional function field
	of characteristic $0$.
	If the $abcd$ conjecture holds, then there is a constant
	$C(d, K) > 0$ such that for any polynomial $f \in K[z]$ of degree
	$d$ and for all $P \in K$,
	either $\hat{h}_f(P) = 0$ or
	\[
		\hat{h}_f(P) \geq C(d, K) \max\set{1, \hat{h}_\crit(f)}.
	\]
\end{theorem}


\section{Acknowledgments}
The author would like to thank the organizers Gaetan Bisson,
Philippe Lebacque, and Roger Oyono of the \textit{G\'{e}om\'{e}trie
alg\'{e}brique, Th\'{e}orie des nombres et Applications}
conference at the University of French Polynesia in August 2021.
The author also benefited from helpful discussions with Michelle Manes
at the conference and discussions with Chatchawan Panraksa
about his unpublished note mentioned in this article.
The author is also grateful to the anonymous referee for
their comments and suggestions.

This work was supported by the
National Science Foundation Graduate Research
Fellowship Program under Grant No. DGE-1644869.


\bibliography{bibliography}{}
\bibliographystyle{amsalpha}

\end{document}